\documentclass[11pt]{article}
\usepackage{amsfonts,latexsym,rawfonts,amsmath,amssymb,amsthm}
\newcommand{\bthm}[2]{\vskip 8pt\noindent\bf #1\hskip 2pt\bf#2\it \hskip 8pt}
\newcommand{\ethm}{\vskip 8pt\rm}
\numberwithin{equation}{section}
\newtheorem{theorem}{Theorem}[section]
\newtheorem{lem}[theorem]{Lemma}
\newtheorem{thm}[theorem]{Theorem}
\newtheorem{pro}[theorem]{Proposition}
\newtheorem{cor}[theorem]{Corollary}
\def\s{\,\,\,\,}

\def\mv{1.8ex}
\def\dint{\displaystyle{\int}}

\title{A sharp Trudinger-Moser type inequality \\for unbounded domains in $\mathbb{R}^n$}
\author{Yuxiang Li and Bernhard Ruf}
\date{{}}

\begin{document}

\maketitle

\begin{abstract}
The Trudinger-Moser inequality states that for functions $u \in
H_0^{1,n}(\Omega)$ ($\Omega \subset \mathbb R^n$ a bounded domain)
with $\int_\Omega |\nabla u|^ndx  \le 1$ one has $\int_\Omega
(e^{\alpha_n|u|^{\frac n{n-1}}}-1)dx \le c\, |\Omega|$, with $c$
independent of $u$. Recently, the second author has shown that for
$n = 2$ the bound $c\, |\Omega| $ may be replaced by a uniform
constant $d$ independent of $\Omega$ if the Dirichlet norm is
replaced by the Sobolev norm, i.e. requiring $\int_\Omega (|\nabla
u|^n + |u|^n )dx \le 1$. We extend here this result to arbitrary
dimensions $n > 2$. Also, we prove that for $\Omega = \mathbb R^n$
the supremum of $\int_{\mathbb R^n} (e^{\alpha_n|u|^{\frac
n{n-1}}}-1)dx$ over all such functions is attained. The proof is
based on a blow-up procedure.\
\end{abstract}
\par \bigskip

\section{Introduction}

Let $H^{1,p}_0(\Omega)$, $\Omega\subseteq\mathbb{R}^n$, be the usual
Sobolev space, i.e. the completion of $C_0^\infty(\Omega)$ with the
norm
$$\|u\|_{H^{1,p}(\Omega)}
=\left(\dint_\Omega(|\nabla u|^p+|u|^p)dx\right)^\frac{1}{p}.$$ It
is well-known that
$$\begin{array}{ll}
        H^{1,p}_0(\Omega)\subset L^\frac{pn}{n-p}(\Omega)&if\s 1\leq p<n\\[1.5ex]
        H^{1,p}_0(\Omega)\subset L^\infty(\Omega)&if\s n<p
  \end{array}$$
The case $p=n$ is the limit case of these embeddings and it is known
that
$$H^{1,n}_0(\Omega)\subset L^q(\Omega)\s for\s n\leq q < +\infty.$$

When $\Omega$ is a bounded domain, we usually use the Dirichlet norm
$\|u\|_D=(\int|\nabla u|^ndx)^\frac{1}{n}$ in  place of  $\|\cdot
\|_{H^{1,n}}$. In this case, we have the famous Trudinger-Moser
inequality (see \cite{P}, \cite{T}, \cite{M}) for the limit case
$p=n$ which states that
\begin{equation}\label{1.1}
   \sup_{\|u\|_D\leq 1}\dint_\Omega(e^{\alpha|u|^\frac{n}{n-1}}-1)dx=c(\Omega,\alpha)
    \left\{\begin{array}{ll}
           <+\infty& when\s \alpha\leq\alpha_n\\[1.5ex]
            =+\infty& when\s \alpha>\alpha_n
        \end{array}\right.
\end{equation}
where $\alpha_n=n\omega_{n-1}^\frac{1}{n-1}$, and $\omega_{n-1}$ is
the measure of the unit sphere in $\mathbb{R}^n$. The
Trudinger-Moser result has been extended to Sobolev spaces of higher
order and Soboleve spaces over compact manifolds (see \cite{A},
\cite{Fo}). Moreover, for any bounded $\Omega$, the constant
$c(\Omega,\alpha_n)$ can be attained. For the attainability, we
refer to \cite{C-C}, \cite{F}, \cite{Lin}, \cite{L}, \cite{L2},
\cite{d-d-R}, \cite{L3}.

Another interesting extension of (\ref{1.1}) is to construct
Trudinger-Moser type inequalities on unbounded domains. When $n=2$,
this has been done by B. Ruf in \cite{R}. On the other hand, for an
unbounded domain in $\mathbb{R}^n$, S. Adachi and K. Tanaka
(\cite{A-T}) get a weaker result. Let
$$\Phi(t)=e^{t}-\sum\limits_{j=1}^{n-2}
\frac{t^j}{j!}.$$ S. Adachi and K. Tanaka's result says that:

\bthm{Theorem}{A} For any $\alpha\in(0,\alpha_n)$ there is a
constant $C(\alpha)$ such that
\begin{equation}\label{1.2}
\dint_{\mathbb{R}^n}\Phi(\alpha(\frac{|u|}{\|\nabla
u\|_{L^n(\mathbb{R}^n)}})^\frac{n}{n-1})dx \leq
C(\alpha)\frac{\|u\|_{L^n(\mathbb{R}^n)}^n}{\|\nabla
u\|_{L^n(\mathbb{R}^n)}^n}\s \ , \ \hbox{ for} \s u\in
H^{1,n}(\mathbb{R}^n)\setminus\{0\}.
\end{equation}\ethm

In this paper, we shall discuss the critical case $\alpha =
\alpha_n$. More precisely, we prove the following:

\begin{thm}\label{t1.1} There exists a constant $d>0$, s.t. for any domain
$\Omega\subset\mathbb{R}^n$,
\begin{equation}\label{1.3}\sup_{u\in H^{1,n}(\Omega),
\|u\|_{H^{1,n}(\Omega)}\leq
1}\dint_\Omega\Phi(\alpha_n|u|^\frac{n}{n-1})dx\leq d.
\end{equation}
The inequality is sharp: for any $\alpha>\alpha_n$, the supremum is $+\infty$.\\
\end{thm}

We set
$$S=\sup_{u\in H^{1,n}(\mathbb{R}^n),
\|u\|_{H^{1,n}(\mathbb{R}^n)}\leq
1}\dint_{\mathbb{R}^n}\Phi(\alpha_n|u|^\frac{n}{n-1})dx.$$ Further,
we will prove

\begin{thm}\label{t1.2} $S$ is attained. In other words, we can find a function $u\in H^{1,n}
(\mathbb{R}^n)$, with $\|u\|_{H^{1,n}(\mathbb{R}^n)}=1$, s.t.
$$S=\dint_{\mathbb{R}^n}\Phi(\alpha_n|u|^{\frac n{n-1}})dx\ .
$$
\end{thm}

The second part of Theorem \ref{t1.1} is trivial: Given any fixed
$\alpha>\alpha_n$, we take $\beta\in (\alpha_n,\alpha)$. By
(\ref{1.1}) we can find a positive sequence $\{u_k\}$ in
$$\{u\in H^{1,n}_0(B_1):\dint_{B_1}|\nabla u|^ndx=1\},$$
such that
$$\lim_{k\rightarrow+\infty}\dint_{B_1}e^{\beta u_k^\frac{n}{n-1}}=+\infty.$$
By  Lion's Lemma, we get $u_k\rightharpoondown 0$. Then by the
compact embedding theorem, we may assume
$\|u_k\|_{L^p(B_1)}\rightarrow 0$ for any $p>1$. Then,
$\int_{\mathbb{R}^n}(|\nabla u_k|^n+|u_k|^n)dx\rightarrow 1$, and
$$\alpha(\frac{u_k}{\|u_k\|_{H^{1,n}}})^\frac{n}{n-1}>\beta u_k^\frac{n}{n-1}$$
when $k$ are sufficiently large. So, we get
$$\lim_{k\rightarrow+\infty}\dint_{\mathbb{R}^n}
\Phi(\alpha(\frac{u_k}{\|u_k\|_{H^{1,n}}})^\frac{n}{n-1})dx\geq
\lim_{k\rightarrow+\infty}\dint_{B_1}(e^{\beta u_k^\frac{n}{n-1}}-1)dx=+\infty.$$\\

The first part of  Theorem \ref{t1.1} and Theorem \ref{t1.2} will be
proved by blow up analysis. We will use the ideas from \cite{L} and
\cite{L2}  (see also \cite{A-M} and \cite{A-D}). However, in the
unbounded case we do not obtain the strong convergence of $u_k$ in
$L^n(\mathbb{R}^n)$, and so we need more techniques.

Concretely, we will find positive and symmetric functions $u_k\in
H_0^{1,n}(B_{R_k})$  which satisfy
$$\dint_{B_{R_k}}(|\nabla u_k|^n+|u_k|^n)dx=1$$
and
$$\dint_{B_{R_k}}\Phi(\beta_ku_k^\frac{n}{n-1})dx=
\sup_{\int_{B_{R_k}}(|\nabla v|^n+|v|^n) = 1,\ v\in
H^{1,n}_0(B_{R_k})}
  \dint_{B_{R_k}}\Phi(\beta_k|v|^\frac{n}{n-1})dx.$$
Here, $\beta_k$ is an increasing sequence tending to $\alpha_n$, and
$R_k$ is an increasing sequence tending to $+\infty$.

Furthermore, $u_k$ satisfies the following equation:
$$-div|\nabla u_k|^{n-2}\nabla u_k+u_k^{n-1}=\frac{u_k^\frac{1}{n-1}
\Phi'(\beta_ku_k^\frac{n}{n-1})}{\lambda_k},$$ where $\lambda_k$ is
a Lagrange multiplier.

Then, there are two possibilities. If $c_k=\max u_k$ is bounded from
above, then it is easy to see that
$$\lim_{k\rightarrow+\infty}
\dint_{\mathbb{R}^n}(\Phi(\beta_ku_k^\frac{n}{n-1})-
\frac{\beta_k^{n-1}u_k^n}{(n-1)!})dx=\dint_{\mathbb{R}^n}
(\Phi(\alpha_nu^\frac{n}{n-1})-
\frac{\alpha_n^{n-1}u^n}{(n-1)!})dx$$ where $u$ is the weak limit of
$u_k$. It then follows that either $\int_{\mathbb{R}^n}\Phi(\beta_k
u_k^\frac{n}{n-1})dx$ converges to $\int_{\mathbb{R}^n}\Phi(\alpha_n
u^\frac{n}{n-1})dx$, or
$$S\leq\frac{\alpha_n^{n-1}}{(n-1)!}.$$

If $c_k$ is not bounded, the key point of the proof is to show that
$$
\frac{n}{n-1}\beta_kc_k^\frac{1}{n-1}(u_k(r_kx)-c_k)\rightarrow
-n\log(1+c_nr^\frac{n}{n-1})\ ,
$$
locally for a suitably chosen sequence $r_k$ (and with $c_n =
(\frac{\omega_{n-1}}{n})^\frac{1}{n-1}$), and that
$$
c_k^\frac{1}{n-1}u_k\rightarrow G\ ,
$$
on any $\Omega\subset\subset\mathbb{R}^n\setminus\{0\}$, where $G$
is some  Green function.
This will be done in section 3.\\

Then, we will get  in section 4 the following

\begin{pro}\label{t1.3} If
$S$ can not be attained, then
$$S\leq \min\{\frac{\alpha_n^{n-1}}{(n-1)!},\frac{\omega_{n-1}}{n}e^{\alpha_nA+
1+1/2+\cdots+1/(n-1)}\},$$ where $A=\lim\limits_{r\rightarrow
0}(G(r) + \frac{1}{\alpha_n}\log{r^n})$.
\end{pro}

So, to prove the attainability, we only need to show that
$$S>\min\{\frac{\alpha_n^{n-1}}{(n-1)!},\frac{\omega_{n-1}}{n}e^{\alpha_nA+
1+1/2+\cdots+1/(n-1)}\}.$$ In  section 5, we will construct a
function sequence $u_\epsilon$ such that
$$\dint_{\mathbb{R}^n}\Phi(\alpha_nu_\epsilon^\frac{n}{n-1})dx>\frac{\omega_{n-1}}{n}
e^{\alpha_nA+ 1+1/2+\cdots+1/(n-1)}$$ when $\epsilon$ is
sufficiently small. And in the last section we will construct, for
each $n > 2$, a function sequence $u_\epsilon$ such that for
$\epsilon$ sufficiently small
$$\dint_{\mathbb{R}^n}\Phi(\alpha_nu_\epsilon^\frac{n}{n-1})dx>\frac{\alpha_n^{n-1}}{(n-1)!}.$$

Thus, together with Ruf's result of attainability  in \cite{R} for
the case $n=2$, we will get Theorem \ref{t1.2}.

\section{The maximizing sequence}
Let $\{R_k\}$ be an increasing sequence which diverges to infinity,
and $\{\beta_k\}$ an increasing sequence which converges to
$\alpha_n$. By compactness, we can find positive functions $u_k\in
H_0^{1,n}(B_{R_k})$ with $\int_{B_{R_k}}(|\nabla u_k|^n+u_k^n)dx=1$
 such that
$$\dint_{B_{R_k}}\Phi(\beta_ku_k^\frac{n}{n-1})dx=
\sup_{\int_{B_{R_k}}(|\nabla v|^n+|v|^n) = 1,\ v\in
H^{1,n}_0(B_{R_k})}
  \dint_{B_{R_k}}\Phi(\beta_k|v|^\frac{n}{n-1})dx.$$
Moreover, we may assume that $\int_{\mathbb{R}^n}\Phi(\beta_k
u_k^\frac{n}{n-1})dx=\int_{B_{R_k}}\Phi(\beta_ku_k^\frac{n}{n-1})dx$
is  increasing.

\begin{lem}\label{t2.1}
Let $u_k$ as above. Then

a) $u_k$ is a maximizing sequence for $S$;

b) $u_k$ may be chosen to be radially symmetric and decreasing.
\end{lem}

\proof a) Let $\eta$ be a cut-off function which is 1 on $B_1$ and 0
on $\mathbb{R}^n\setminus B_2$. Then given any $\varphi\in
H^{1,n}(\mathbb{R}^n)$ with $\int_{\mathbb{R}^n}(|\nabla
\varphi|^n+|\varphi|^n)dx=1$, we have
$$\tau^n(L):=\dint_{\mathbb{R}^n}(|\nabla\eta(\frac{x}{L})\varphi|^n
+|\eta(\frac{x}{L})\varphi|^n)dx\rightarrow 1,\s as\s
L\rightarrow+\infty.$$ Hence for a fixed $L$ and $R_k > 2L$
$$
\dint_{B_L}\Phi(\beta_k|\frac{\varphi}{\tau(L)}|^\frac{n}{n-1})dx
\leq\dint_{B_{2L}}\Phi(\beta_k|\frac{\eta(\frac{x}{L})\varphi}{\tau(L)}|^\frac{n}{n-1})dx
\leq\dint_{B_{R_k}} \Phi(\beta_ku_k^\frac{n}{n-1})dx
$$
By the Levi Lemma, we then have
$$\dint_{B_L}\Phi(\alpha_n|\frac{\varphi}{\tau(L)}|^\frac{n}{n-1})dx
\leq\lim_{k\rightarrow+\infty}\dint_{\mathbb{R}^n}
\Phi(\beta_ku_k^\frac{n}{n-1})dx.$$ Then, letting
$L\rightarrow+\infty$, we get
$$\dint_{\mathbb{R}^n}\Phi(\alpha_n|\varphi|^\frac{n}{n-1})dx
\leq \lim_{k\rightarrow+\infty}\dint_{\mathbb{R}^n}
\Phi(\beta_ku_k^\frac{n}{n-1})dx.$$ Hence, we get
$$\lim_{k\rightarrow+\infty}\dint_{\mathbb{R}^n}
\Phi(\beta_ku_k^\frac{n}{n-1})dx =\sup_{\int_{\mathbb{R}^n} (|\nabla
v|^n+|v|^n) = 1,\ v\in H^{1,n}(\mathbb{R}^n)}
  \dint_{\mathbb{R}^n}\Phi(\alpha_n|v|^\frac{n}{n-1})dx.$$

b) Let $u_k^*$ be the radial rearrangement of $u_k$, then we have
$$\tau_k^n:=\dint_{B_{R_k}}(|\nabla u_k^*|^n+{u_k^*}^n)dx\leq
\dint_{B_{R_k}}(|\nabla u_k|^n+u_k^n)dx=1.$$ It is well-known that
$\tau_k=1$ iff $u_k$ is radial. Since
$$\dint_{B_{R_k}}\Phi(\beta_k{u_k^*}^\frac{n}{n-1})dx
=\dint_{B_{R_k}}\Phi(\beta_ku_k^\frac{n}{n-1})dx,$$ we have
$$\dint_{B_{R_k}}\Phi(\beta_k(\frac{u_k^*}{\tau_k})^\frac{n}{n-1})dx
\geq \dint_{B_{R_k}}\Phi(\beta_ku_k^\frac{n}{n-1})dx,$$ and "$=$\, "
holds iff $\tau_k=1$. Hence $\tau_k=1$ and
$$\dint_{B_{R_k}}\Phi(\beta_k{u_k^*}^\frac{n}{n-1})dx
=\sup_{\int_{B_{R_k}}(|\nabla v|^n+|v|^n) = 1,\ v\in
H^{1,n}_0(B_{R_k})}
  \dint_{B_{R_k}}\Phi(\beta_k|v|^\frac{n}{n-1})dx.$$
So, we can assume $u_k=u_k(|x|)$, and $u_k(r)$ is decreasing.

$\hfill\Box$\\

Assume now $u_k\rightharpoondown u$. Then, to prove Theorem 1.1 and
1.2, we only need to show that
$$\lim_{k\rightarrow+\infty}\dint_{\mathbb{R}^n}\Phi(\beta_k
u_k^\frac{n}{n-1})dx
=\dint_{\mathbb{R}^n}\Phi(\alpha_nu^\frac{n}{n-1})dx.$$

\section{Blow up analysis}

By the definition of $u_k$ we have the equation
\begin{equation}\label{3.1}
-div|\nabla u_k|^{n-2}\nabla u_k+u_k^{n-1}=\frac{u_k^\frac{1}{n-1}
\Phi'(\beta_ku_k^\frac{n}{n-1})} {\lambda_k}\ ,
\end{equation}
where $\lambda_k$ is the constant satisfying
$$\lambda_k=\dint_{B_{R_k}}u_k^{\frac{n}{n-1}}\Phi'(\beta_ku_k^\frac{n}{n-1})dx.$$

First, we need to prove the following:

\begin{lem}\label{t3.1} $\inf\limits_k\lambda_k>0$.\end{lem}

\proof Assume $\lambda_k\rightarrow 0$. Then
$$
\int_{\mathbb{R}^n}u_k^ndx\leq
C\int_{\mathbb{R}^n}u_k^\frac{n}{n-1}\Phi'(\beta_ku_k^\frac{n}{n-1})dx\leq
C\lambda_k\rightarrow 0\ .
$$
Since $u_k(|x|)$ is decreasing, we have $u_k^n(L)|B_L|\leq
\int_{B_L}u_k^n\leq 1$, and then
\begin{equation}\label{3.2a}
u_k(L)\leq\frac{n}{\omega_nL^n}.
\end{equation}
Set $\epsilon=\frac{n}{\omega_nL^n}$. Then $u_k(x)\leq\epsilon$ for
any $x\notin B_L$, and hence we have, using the form of $\Phi$, that
$$\dint_{\mathbb{R}^n\setminus B_L}\Phi(\beta_ku_k^\frac{n}{n-1})dx\leq
C\dint_{\mathbb{R}^n\setminus B_L}u_k^ndx\leq C\lambda_k \rightarrow
0\ .
$$
And on $B_L$, since $u_k\rightarrow 0$ in $L^q(B_L)$ for any $q>1$,
we have by Lebesgue
$$\begin{array}{ll}
     \lim\limits_{k\rightarrow+\infty}\dint_{B_L}\Phi(\beta_ku_k^\frac{n}{n-1})dx
       &\leq\lim\limits_{k\rightarrow+\infty}\left[\dint_{B_L}Cu_k^\frac{n}{n-1}\Phi'(\beta_k
       u_k^\frac{n}{n-1})dx
         +\dint_{\{x\in B_L: u_k(x)\leq 1\}}\Phi(\beta_ku_k^\frac{n}{n-1})dx\right]\\[\mv]
       &\leq\lim\limits_{k\rightarrow+\infty}C\lambda_k+
         \dint_{B_L}\Phi(0)dx\\[\mv]
       &=0.
  \end{array}$$
This is impossible.

$\hfill\Box$\\

We denote $c_k=\max u_k=u_k(0)$. Then we have

\begin{lem}\label{t3.2}  If \ $\sup\limits_kc_k<+\infty$, then

i) \ Theorem 1.1 holds; \vspace{0.2cm}

ii) if S is not attained, then
$$S\leq \frac{\alpha_n^{n-1}}{(n-1)!}\ .
$$
\end{lem}

\proof If $\sup_kc_k<+\infty$, then $u_k\rightarrow u$ in
$C^1_{loc}(\mathbb{R}^n)$. By (\ref{3.2a}), we are able to find $L$
s.t. $u_k(x)\leq\epsilon$ for $x\notin{B_L}$. Then
$$\dint_{\mathbb{R}^n\setminus B_L}(\Phi(\beta_ku_k^\frac{n}{n-1})-\frac{\beta_k^{n-1}
u_k^n}{(n-1)!})dx \leq C\dint_{\mathbb{R}^n\setminus
B_L}u_k^\frac{n^2}{n-1}dx\leq C\epsilon^{\frac{n^2}{n-1}
-n}\dint_{\mathbb{R}^n}u_k^ndx\leq C\epsilon^{\frac{n^2}{n-1} -n}.$$
Letting $\epsilon\rightarrow 0$, we get
$$\lim_{k\rightarrow+\infty}
\dint_{\mathbb{R}^n}(\Phi(\beta_ku_k^\frac{n}{n-1})-\frac{\beta_k^{n-1}u_k^n}{(n-1)!})dx
=\dint_{\mathbb{R}^n}(\Phi(\alpha_nu^\frac{n}{n-1})-\frac{\alpha_n^{n-1}u^n}{(n-1)!})dx.$$
Hence
\begin{equation}\label{3.2}
\lim_{k\rightarrow+\infty}\dint_{\mathbb{R}^n}\Phi(\beta_ku_k^\frac{n}{n-1})=
\dint_{\mathbb{R}^n}\Phi(\alpha_nu^\frac{n}{n-1})dx+\frac{\alpha_n^{n-1}}{(n-1)!}
\lim_{k\rightarrow+\infty} \dint_{\mathbb{R}^n}(u_k^n -u^n)dx.
\end{equation}

When $u=0$, we can deduce from (\ref{3.2}) that
$$S\leq\frac{\alpha_n^{n-1}}{(n-1)!}.$$

Now, we assume $u\neq 0$. Set
$$\tau^n=\lim_{k\rightarrow+\infty}
\frac{\dint_{\mathbb{R}^n}u_k^ndx} {\dint_{\mathbb{R}^n}u^ndx}.$$ By
the Levi Lemma, we have $\tau\geq 1$.

Let $\tilde{u}=u(\frac{x}{\tau})$. Then, we have
$$\dint_{\mathbb{R}^n}|\nabla \tilde{u}|^ndx=\dint_{\mathbb{R}^n}
|\nabla
u|^ndx\leq\lim_{k\rightarrow+\infty}\dint_{\mathbb{R}^n}|\nabla
u_k|^ndx,$$ and
$$\dint_{\mathbb{R}^n}\tilde{u}^ndx=\tau^n\dint_{\mathbb{R}^n}u^ndx=
\lim_{k\rightarrow+\infty}\dint_{\mathbb{R}^n}u_k^ndx.$$ Then
$$\dint_{\mathbb{R}^n}(|\nabla \tilde{u}|^n+\tilde{u}^n)dx\leq
\lim_{k\rightarrow+\infty}\dint_{\mathbb{R}^n}(|\nabla u_k|^n
+u_k^n)dx=1.$$ Hence, we have by (\ref{3.2})
$$\begin{array}{lll}
    S&\geq&\dint_{\mathbb{R}^n}\Phi(\alpha_n\tilde{u}^\frac{n}{n-1})dx\\[\mv]
      &=&\tau^n\dint_{\mathbb{R}^n}\Phi(\alpha_nu^\frac{n}{n-1})dx\\[\mv]
      &=&\left[\dint_{\mathbb{R}^n}\Phi(\alpha_nu^\frac{n}{n-1})dx+
          (\tau^n-1)\dint_{\mathbb{R}^n}\frac{\alpha_n^{n-1}}{(n-1)!}u^ndx\right]
      +(\tau^n-1)\dint_{\mathbb{R}^n}(\Phi(\alpha_nu^\frac{n}{n-1})
           -\frac{\alpha_n^{n-1}}{(n-1)!}u^n)dx\\[\mv]
      &=&\lim\limits_{k\rightarrow+\infty}\dint_{\mathbb{R}^n}\Phi(\beta_ku_k^\frac{n}{n-1})dx
         +(\tau^n-1)\dint_{\mathbb{R}^n}(\Phi(\alpha_nu^\frac{n}{n-1})
           -\frac{\alpha_n^{n-1}}{(n-1)!}u^n)dx\\[\mv]
      &=&S+(\tau^n-1)\dint_{\mathbb{R}^n}(\Phi(\alpha_nu^\frac{n}{n-1})
           -\frac{\alpha_n^{n-1}}{(n-1)!}u^n)dx
  \end{array}$$

Since
$\Phi(\alpha_nu^\frac{n}{n-1})-\frac{\alpha_n^{n-1}}{(n-1)!}u^n>0$,
we have $\tau=1$, and then
$$S=\dint_{\mathbb{R}^n}\Phi(\alpha_nu^\frac{n}{n-1})dx.$$
So, $u$ is an extremal function.

$\hfill\Box$ \vspace{2ex}

From now on, we assume $c_k\rightarrow+\infty$. We perform a blow-up
procedure:

\noindent We define
$$r_k^n=\frac{\lambda_k}{c_k^{\frac{n}{n-1}}e^{\beta_kc_k^\frac{n}{n-1}}}.$$
By (\ref{3.2a}) we can find a sufficiently large $L$ such that
$u_k\leq 1$ on $\mathbb{R}^n\setminus B_L$. Then
$$\dint_{B_L}|\nabla (u_k-u_k(L))^+|^ndx\leq 1$$
and hence, by (\ref{1.1}), we have
$$\dint_{B_L}e^{\alpha_n[(u_k-u_k(L))^+]^\frac{n}{n-1}}\leq C(L).$$
Clearly, for any $p<\alpha_n$ we can find a constant $C(p)$, s.t.
$$pu_k^\frac{n}{n-1}\leq \alpha_n[(u_k-u_k(L))^+]^\frac{n}{n-1}+C(p),$$
and then we get
$$\dint_{B_L}e^{pu_k^\frac{n}{n-1}}dx<C=C(L,p).$$
Hence,
$$\begin{array}{ll}
     \lambda_ke^{-\frac{\beta_k}{2}c_k^\frac{n}{n-1}}
        &= e^{-\frac{\beta_k}{2}c_k^\frac{n}{n-1}}
          \left[\dint_{\mathbb{R}^n\setminus B_L}u_k^\frac{n}{n-1}\Phi'(\beta_k
          u_k^\frac{n}{n-1})
          dx+\dint_{B_L}u_k^\frac{n}{n-1}\Phi'(\beta_ku_k^\frac{n}{n-1})dx\right]\\[\mv]
        &\leq C\dint_{\mathbb{R}^n\setminus B_L}
             u_k^ndx\; e^{-\frac{\beta_k}{2}c_k^\frac{n}{n-1}}
             +\dint_{B_L}e^{\frac{\beta_k}{2}u_k^\frac{n}{n-1}}
             u_k^\frac{n}{n-1}dx.
  \end{array}$$
Since $u_k$ converges strongly in $L^q(B_L)$ for any $q>1$, we get
$\lambda_k\leq Ce^{\frac{\beta_k}{2}c_k^\frac{n}{n-1}}$, and hence
$$r_k^n\leq Ce^{-\frac{\beta_k}{2}c_k^\frac{n}{n-1}}.$$

Now, we set
$$v_k(x)=u_k(r_kx),\s w_k(x)=\frac{n}{n-1}\beta_kc_k^\frac{1}{n-1}
(v_k-c_k),$$ where $v_k$ and $w_k$ are  defined on $\Omega_k=\{x\in
\mathbb{R}^n:r_kx\in B_1\}$. Using the definition of $r_k^n$ and
(\ref{3.1}) we have
$$-div|\nabla w_k|^{n-2}\nabla w_k=\frac{v_k^\frac{1}{n-1}
}{c_k^\frac{1
  }{n-1}}(\frac{n}{n-1}\beta_k)^{n-1}e^{\beta_k(v_k^\frac{n}{n-1}-
  c_k^\frac{n}{n-1})}+O(r_k^nc_k^n).$$

By Theorem 7 in \cite{S}, we know that $osc_{B_R}\omega_k\leq C(R)$
for any $R>0$. Then from the result in  $\cite{T}$ (or $\cite{D}$),
it follows that $\|w_k\|_{C^{1,\delta}(B_R)}<C(R)$. Therefore $w_k$
converges in $C^1_{loc}$ and $v_k-c_k\rightarrow 0$ in $C^1_{loc}$.

Since
$$ v_k^\frac{n}{n-1} = c_k^{\frac{n}{n-1}}(1+\frac{v_k-c_k}{c_k})^{\frac{n}{n-1}}
=c_k^{\frac{n}{n-1}}(1+\frac{n}{n-1}
\frac{v_k-c_k}{c_k}+O(\frac{1}{c_k^2}))\ ,
$$
we get $\beta_k(v_k^\frac{n}{n-1}-c_k^\frac{n}{n-1})\rightarrow w$
in $C^0_{loc}$, and so we have
\begin{equation}\label{3.3}
-div|\nabla w|^{n-2}\nabla w=(\frac{n\alpha_n}{n-1})^{n-1}e^w \ ,
\end{equation}
with $$w(0)=0=\max w.$$

Since $\omega$ is radially symmetric and decreasing, it is easy to
see that (\ref{3.3}) has only one solution. We can check that
$$w(x)=-n\log(1+c_n|x|^\frac{n}{n-1}),\s \hbox{and} \s \dint_{\mathbb{R}^n}e^wdx=1,$$
where $c_n=(\frac{\omega_{n-1}}{n})^\frac{1}{n-1}$. Then,
\begin{equation}\label{3.4}
\lim_{L\rightarrow +\infty}\lim_{k\rightarrow +\infty}
\dint_{B_{Lr_k}}\frac{u_k^\frac{n}{n-1}}{\lambda_k}e^{
\beta_ku_k^\frac{n}{n-1}}dx=\lim_{L\rightarrow
+\infty}\dint_{B_L}e^wdx=1.
\end{equation}
\vspace{2ex}

For $A>1$, let $u_k^A=\min\{u_k,\frac{c_k}{A}\}$. We have

\begin{lem}\label{t3.3} For any $A>1$, there holds
\begin{equation}\label{3.5}
\limsup_{k\rightarrow+\infty}\dint_{\mathbb{R}^n} (|\nabla
u_k^A|^n+|u_k^A|^n)dx\leq\frac{1}{A}.
\end{equation}
\end{lem}

\proof Since $|\{x:u_k\geq\frac{c_k}{A}\}|\,
|\frac{c_k}{A}|^n\leq\int_{\{u_k\geq \frac{c_k}{A}\}}u_k^n\leq 1$,
we can find a sequence $\rho_k\rightarrow 0$ s.t.
$$\{x:u_k\geq\frac{c_k}{A}\}\subset B_{\rho_k}.$$
Since $u_k$ converges in $L^p(B_1)$ for any $p>1$, we have
$$\lim_{k\rightarrow+\infty}\dint_{\{u_k>\frac{c_k}{A}\}}|u_k^A|^pdx
\leq\lim_{k\rightarrow+\infty}\dint_{\{u_k>\frac{c_k}{A}\}}u_k^pdx=
0,$$ and
$$\lim_{k\rightarrow+\infty}\dint_{\mathbb{R}^n}(u_k-\frac{c_k}{A})^+u_k^pdx= 0$$
for any $p>0$.

Hence, testing equation (\ref{3.1}) with $(u_k - \frac{c_k}A)^+$, we
have
$$\begin{array}{ll}
      \dint_{\hspace{-0.15cm}\mathbb{R}^n}\left(|\nabla (u_k-\frac{c_k}{A})^+|^n+(u_k-\frac{c_k}{A})^+
          u_k^{n-1}\right)dx
        &=\dint_{\hspace{-0.15cm}\mathbb{R}^n}
        (u_k-\frac{c_k}{A})^+\frac{u_k^{\frac{1}{n-1}}}{\lambda_k}e^{\beta_ku_k^\frac{n}
          {n-1}}dx+o(1)\vspace{0.2cm}\\[\mv]
      &\geq\dint_{\hspace{-0.15cm} B_{Lr_k}}(u_k-\frac{c_k}{A})^+
        \frac{u_k^\frac{1}{n-1}}{\lambda_k}e^{\beta_ku_k^\frac{n}{n-1}}dx+o(1)\vspace{0.3cm}\\[\mv]
      &=\dint_{\hspace{-0.15cm} B_L}\frac{v_k-c_k/A}{c_k}(\frac{v_k-c_k}{c_k}+1)^\frac{1}{n-1}
        e^{w_k+o(1)}dx+o(1).\\[\mv]
  \end{array}$$
Hence
$$\liminf\limits_{k\rightarrow+\infty}\dint_{\mathbb{R}^n}\left(|\nabla (u_k-\frac{c_k}
{A})^+|^n +(u_k-\frac{c_k}{A})^+u_k^{n-1}\right)dx
\geq\frac{A-1}{A}\int_{B_L}e^wdx.$$ Letting $L\rightarrow+\infty$,
we get
$$\liminf_{k\rightarrow+\infty}\dint_{\mathbb{R}^n}\left(|\nabla (u_k-\frac{c_k}{A})^+|^n
+(u_k-\frac{c_k}{A})^+u_k^{n-1}\right)dx\geq\frac{A-1}{A}.
$$
Now observe that
$$\begin{array}{lll}
   \dint_{\mathbb{R}^n}(|\nabla u_k^A|^n+|u_k^A|^n)dx
     &=&1-\dint_{\mathbb{R}^n}\left(|\nabla(u_k-\frac{c_k}{A})^+|^n+(u_k-\frac{c_k}{A})^+
       u_k^{n-1}\right)dx \vspace{0.2cm}\\[\mv]
     &&+\dint_{\mathbb{R}^n}(u_k-\frac{c_k}{A})^+u_k^{n-1}dx-
       \dint_{\{u_k>\frac{c_k}{A}\}}u_k^ndx+\dint_{\{u_k>\frac{c_k}{A}\}}|u_k^A|^ndx\vspace{0.2cm}\\[\mv]
     &\leq&1-(1-\frac{1}{A})+o(1).
  \end{array}
$$
Hence, we get this Lemma.

$\hfill\Box$\\

\begin{cor}\label{t3.4}
We have
$$\lim_{k\rightarrow+\infty}\dint_{\mathbb{R}^n\setminus B_\delta}(|\nabla u_k|^n
+u_k^n)dx=0\ ,
$$
for any $\delta>0$, and then $u=0$.
\end{cor}
\proof Letting $A\rightarrow+\infty$, then for any constant $c$, we
have
$$\int_{\{u_k\leq c\}}(|\nabla u_k|^n+u_k^n)dx\rightarrow 0.$$
So  we get this Corollary.

$\hfill\Box$\\

\begin{lem}\label{t3.5} We have
\begin{equation}\label{3.7}
\lim_{k\rightarrow+\infty}\dint_{\mathbb{R}^n}\Phi(\beta_ku_k^\frac{n}{n-1})dx
\leq
\lim_{L\rightarrow+\infty}\lim_{k\rightarrow+\infty}\dint_{B_{Lr_k}}
(e^{\beta_ku_k^\frac{n}{n-1}}-1)dx=
\limsup_{k\rightarrow\infty}\frac{\lambda_k}{c^\frac{n}{n-1}_k},
\end{equation}
and consequently
\begin{equation}\label{3.8a}
\frac{\lambda_k}{c_k}\rightarrow +\infty\ ,\s and\s
\sup\limits_{k}\frac{c_k^\frac{n}{n-1}}{\lambda_k}<+\infty\ .
\end{equation}
\end{lem}

\proof  We have
$$\begin{array}{ll}
      \dint_{\mathbb{R}^n}\Phi(\beta_ku_k^\frac{n}{n-1}dx)& \le \dint_{\{u_k\leq
        \frac{c_k}{A}\}}
        \Phi(\beta_ku_k^\frac{n}{n-1})dx+\dint_{\{u_k> \frac{c_k}{A}\}}
        \Phi'(\beta_ku_k^\frac{n}{n-1})dx\\[\mv]
      &\leq \dint_{\mathbb{R}^n}\Phi(\beta_k(u_k^A)^\frac{n}{n-1})dx+A^\frac{n}{n-1}
       \frac{\lambda_k}{c_k^\frac{n}{n-1}}
       \dint_{\mathbb{R}^n}\frac{u_k^\frac{n}{n-1}}{\lambda_k}\Phi'(\beta_ku_k^\frac{n}
       {n-1})dx\ .
  \end{array}$$

Applying (\ref{3.2a}), we can find $L$ such that $u_k\leq 1$ on
$\mathbb{R}^n \setminus B_{L}$. Then by Corollary \ref{t3.4} and the
form of $\Phi$, we have
\begin{equation}\label{3.8b}
\lim_{k\rightarrow+\infty}\dint_{\mathbb{R}^n\setminus B_L}
\Phi(p\beta_k(u_k^A)^\frac{n}{n-1})dx \le \lim_{k\to\infty} C(p)
\dint_{\mathbb{R}^n\setminus B_L}u_k^ndx = 0\,
\end{equation}
for any $p>0$.

Since by Lemma \ref{t3.3}\ $\limsup\limits_{k\rightarrow+\infty}
\int_{\mathbb{R}^n}(|\nabla u_k^A|^n+|u_k^A|^n)dx \leq\frac{1}{A}<1$ when
$A>1$, it follows from (\ref{1.1}) that
$$\sup_k\dint_{B_L}
e^{p'\beta_k((u_k^A-u_k(L))^+)^\frac{n}{n-1}}dx<+\infty$$
for any $p'<A^\frac{1}{n-1}$. Since for any $p<p'$
$$p(u_k^A)^\frac{n}{n-1}\leq
p'((u_k^A-u_k(L))^+)^\frac{n}{n-1}+C(p,p'),$$
 we have
\begin{equation}\label{p}
\sup_k\dint_{B_L}\Phi(p\beta_k(u_k^A)^{\frac{n}{n-1}})dx<+\infty
\end{equation}
for any $p<A^\frac{1}{n-1}$. Then on $B_L$, by the weak compactness
of Banach space, we get
$$\lim_{k\rightarrow+\infty}\dint_{B_L}\Phi(\beta_k(u_k^A)^\frac{n}{n-1})dx=
\dint_{B_L}\Phi(0)dx=0\ .
$$

Hence we have
$$\begin{array}{ll}
     \lim\limits_{k\rightarrow+\infty}\dint_{\mathbb{R}^n}\Phi(\beta_ku_k^\frac{n}{n-1})dx&
     \displaystyle \leq
           \lim\limits_{L \to \infty} \lim\limits_{k\rightarrow+\infty}A^\frac{n}{n-1}\frac{\lambda_k}
               {c_k^\frac{n}{n-1}}\dint_{B_{L}}\frac{u_k^\frac{n}{n-1}}{\lambda_k}
                \Phi'(\beta_ku_k^\frac{n}{n-1})dx+C\epsilon \vspace{0.2cm}\\[\mv]
         &= \displaystyle \lim\limits_{k\rightarrow+\infty}A^\frac{n}{n-1}\frac{\lambda_k}
               {c_k^\frac{n}{n-1}}+C\epsilon.
  \end{array}$$
As $A\rightarrow 1$ and $\epsilon\rightarrow 0$ we obtain
(\ref{3.7}).

If $\frac{\lambda_k}{c_k}$ was bounded or
$\sup\limits_{k}\frac{c_k^\frac{n}{n-1}}{\lambda_k} =+\infty$, it
would follow from (\ref{3.7}) that
$$\sup\limits_{\int_{\mathbb{R}^n}(|\nabla
v|^n+|v|^n)dx=1,v\in H^{1,n}(\mathbb{R}^n)}
\dint_{\mathbb{R}^n}\Phi(\alpha_n|v|^\frac{n}{n-1})dx=0,$$ which is
impossible.

$\hfill\Box$\\

\begin{lem}\label{t3.6} We have that $c_k\frac{u_k^\frac{1}{n-1}}{
  \lambda_k}\Phi'(\beta_ku_k^\frac{n}{n-1})$ converges to $\delta_0$ weakly, i.e.
for any $\varphi\in\mathcal{D}(\mathbb{R}^n)$ we have
$$\lim_{k\rightarrow+\infty}\dint_{\mathbb{R}^n}\varphi \; c_k\frac{u_k^\frac{1}{n-1}}{
  \lambda_k}\Phi'(\beta_ku_k^\frac{n}{n-1})dx=\varphi(0).$$
\end{lem}

\proof Suppose $supp\; \varphi\subset B_{\rho}$. We split the
integral
$$\begin{array}{ll}
    \dint_{B_\rho}\varphi\; \frac{c_ku_k^\frac{1}{n-1}}{\lambda_k}\Phi'(\beta_ku_k^\frac{n}{n-1})dx
    & \le
       \dint_{\{u_k\geq \frac{c_k}{A}\}\setminus B_{Lr_k}}\cdots+
       \dint_{B_{Lr_k}}\cdots+\dint_{\{u_k<\frac{c_k}{A}\}}\cdots \vspace{0.4cm}\\[\mv]
    &=I_1+I_2+I_3\ .
   \end{array}$$
We have
$$I_1\leq A \|\varphi\|_{C^0}\dint_{\mathbb{R}^n\setminus B_{Lr_k}}\frac{
  u_k^\frac{n}{n-1}}{\lambda_k}\Phi'(\beta_ku_k^\frac{n}{n-1})dx=
  A\|\varphi\|_{C^0}(1-\dint_{B_L}e^{w_k+o(1)}dx)\ ,$$
and
$$I_2=\dint_{B_L}\varphi(r_kx)\frac{c_k(c_k+(v_k-c_k))^\frac{1}{n-1}}
  {c_k^\frac{n}{n-1}}e^{w_k+o(1)}dx
     =\varphi(0)\dint_{B_L}e^wdx+o(1)=\varphi(0)+o(1)\ .$$

By (\ref{3.8b}) and (\ref{p}) we have
$$\dint_{\mathbb{R}^n}\Phi(p\beta_k|u_k^A|^\frac{n}{n-1})dx<C$$
for any $p<A^\frac{1}{n-1}$. We set $\frac{1}{q}+\frac{1}{p}=1$.
Then we get by (\ref{3.8a})
$$
I_3 = \dint_{\{u_k\leq\frac{c_k}{A}\}}\varphi\;
c_k\frac{u_k^\frac{1}{n-1}}{
  \lambda_k}\Phi'(\beta_ku_k^\frac{n}{n-1})dx\leq
\frac{c_k}{\lambda_k}\|\varphi \|_{C^0}\
\|u_k^\frac{1}{n-1}\|_{L^q(\mathbb R^n)}\
 \|e^{\beta_k|u_k^A|^\frac{n}{n-1}} \|_{L^p(\mathbb R^n)}\rightarrow 0.$$

Letting $L\rightarrow+\infty$, we deduce now that
$$\lim_{k\rightarrow+\infty}\dint_{\mathbb{R}^n}\varphi\; \frac{c_ku_k^\frac{1}{n-1}}
{\lambda_k}\Phi'(\beta_ku_k^\frac{n}{n-1})dx=\varphi(0).$$

$\hfill\Box$\\

\begin{pro}\label{t3.7} On any $\Omega\subset\subset
\mathbb{R}^n\setminus\{0\}$, we have that $c_k^\frac{1}{n-1}u_k$
converges to $G$ in $C^1(\Omega)$, where $G\in
C^{1,\alpha}_{loc}(\mathbb{R}^n\setminus\{0\})$ satisfies the
following equation:
\begin{equation}\label{3.9}
-div|\nabla G|^{n-2}\nabla G+G^{n-1}=\delta_0 \ .
\end{equation}
\end{pro}

\proof We set $U_k=c_k^\frac{1}{n-1}u_k$, which satisfy by
(\ref{3.1}) the equations:
\begin{equation}\label{3.10}
-div|\nabla U_k|^{n-2}\nabla
U_k+U_k^{n-1}=\frac{c_ku_k^\frac{1}{n-1}}{\lambda_k}
\Phi'(\beta_ku_k^\frac{n}{n-1})\ .
\end{equation}

For our purpose, we need to prove that
$$
\dint_{B_R}|U_k|^q dx \leq C(q,R)\ ,
$$
where $C(q,R)$ does not depend on $k$. We use the idea in \cite{St}
to prove this statement.

Set $\Omega_t=\{0\leq U_k\leq t\}$, $U_k^t=\min\{U_k ,t\}$. Then we
have
$$\dint_{\Omega_t}(|\nabla U^t_k|^n+|U_k^t|^n)dx\leq
\dint_{\mathbb{R}^n}(-U_k^t\Delta_nU_k+U_k^tU_k^{n-1})=
\dint_{\mathbb{R}^n}U_k^t\frac{c_ku_k^\frac{1}{n-1}}{\lambda_k}\Phi'(\beta_k
u_k^\frac{n}{n-1})dx \leq 2t.$$ Let $\eta$ be a radially symmetric
cut-off function which is 1 on $B_R$ and $0$ on $B_{2R}^c$. Then,
$$\dint_{B_{2R}}|\nabla \eta U_k^t|^ndx\leq C_1(R)+C_2(R)t.$$
Then, when $t$ is bigger than $\frac{C_1(R)}{C_2(R)}$, we have
$$\dint_{B_{2R}}|\nabla\eta U_k^t|^ndx\leq 2C_2(R)\, t\ .
$$
Set $\rho$ such that $U_k(\rho)=t$. Then we have
$$
\inf\left\{\dint_{B_{2R}}|\nabla v|^ndx: v\in H^{1,n}_0(B_{2R})
 \hbox{ and } v|_{B_\rho}=t\right\}\leq 2C_2(R)\, t\ .
 $$
On the other hand, the inf is achieved by
$-t\log{\frac{|x|}{2R}}/\log{ \frac{2R}{\rho}}$. By a direct
computation, we have
$$\frac{\omega_{n-1}t^{n-1}}{(\log\frac{2R}{\rho})^{n-1}}\leq 2C_2(R)\ ,$$
and hence for any $t> \frac{C_1(R)}{C_2(R)}$
$$
|\{x\in B_{2R}:U_k\geq t \}|
 =|B_\rho|\leq C_3(R)e^{-A(R)t}\, ,
$$
where $A(R)$ is a constant only depending on $R$. Then, for any
$\delta<A\, ,$
$$
\dint_{B_R}e^{\delta U_k}dx\leq \sum_{m=0}^\infty\mu(\{m\leq U_k
\leq m+1\})e^{\delta(m+1)}\leq \sum_{m=0}^\infty
e^{-(A-\delta)m}e^\delta\leq C\, .
$$

Then, testing the equation (\ref{3.10}) with the function
$\log{\frac{1+2(U_k-U_k(R))^+}{1+(U_k-U_k(R))^+}}\, ,$ we get
$$
\begin{array}{l}
 \dint_{B_R}\frac{|\nabla U_k|^n}{(1+U_k-U_k(R))(1+2U_k-2U_k(R))}dx\vspace{0.2cm}\\[\mv]
  \s\s\s\s\s\leq \log2\dint_{B_R}\frac{c_ku_k^\frac{1}{n-1}}{\lambda_k}
\Phi'(\beta_ku_k^\frac{n}{n-1})dx-\dint_{B_R}U_k^{n-1}
  \log{\frac{1+2(U_k-U_k(R))}{1+(U_k-U_k(R))}}dx\leq C\, .
\end{array}
$$
Given $q<n$, by Young's Inequality, we have
$$\begin{array}{ll}
   \dint_{B_R}|\nabla U_k|^qdx&\leq\dint_{B_R}\left[\frac{|\nabla U_k|^n}
      {(1+U_k-U_k(R))(1+2U_k-2U_k(R))}+((1+U_k)(
            1+2U_k))^\frac{n}{n-q}\right]dx\vspace{0.3cm}\\[\mv]
            &\leq \dint_{B_R}\left[\frac{|\nabla U_k|^n}{(1+U_k-U_k(R))(1+2U_k-2U_k(R))}
             +Ce^{\delta U_k}\right]dx\, .
            \end{array}$$
Hence, we are able to assume that $U_k$ converges to a function $G$
weakly in $H^{1,p}(B_R)$ for any $R$ and $p<n$. Applying Lemma
\ref{t3.6}, we get (\ref{3.9}).

Hence $U_k$ is bounded in $L^q(\Omega)$ for any $q>0$. By Corollary
\ref{t3.4} and Theorem A, $e^{\beta_ku_k^\frac{n}{n-1}}$ is also
bounded in $L^q(\Omega)$ for any $q>0$. Then, applying Theorem 2.8
in \cite{S}, and the main result in \cite{T} (or \cite{D}), we get
$\|U_k\|_{C^{1,\alpha}(\Omega)}\leq C$. So, $U_k$ converges to $G$
in $C^1(\Omega)$.

$\hfill\Box$\\

For the Green function G we have the following results:
\par \medskip
\begin{lem}\label{t3.8} $G\in C^{1,\alpha}_{loc}(\mathbb{R}^n\setminus\{0\})$ and near
$0$ we can write
\begin{equation}\label{3.11}
G=-\frac{1}{\alpha_n}\log r^n+A+O(r^n\log^n r)\ ;
\end{equation}
here, $A$ is a constant. Moreover, for any $\delta>0$, we have
$$\begin{array}{lll}
    \lim\limits_{k\rightarrow+\infty}\dint_{\mathbb{R}^n\setminus B_\delta}(|\nabla
    c_k^\frac{1}{n-1}u_k|^n+(c_k^\frac{1}{n-1}u_k)^n)dx
    &=&\dint_{\mathbb{R}^n\setminus B_\delta}(|\nabla G|^n+|G|^n)dx\\[\mv]
    &=&G(\delta)(1-\dint_{B_\delta}G^{n-1}dx).
  \end{array}$$
\end{lem}

\proof Slightly modifying the proof in \cite{K-L}, we can prove
$$G=-\frac{1}{\alpha_n}\log r^n+A+o(1).$$
One can refer to \cite{L2} for details. Further, testing the
equation (\ref{3.10}) with 1, we get
$$\omega_{n-1}{G'(r)}^{n-1}r^{n-1}=\dint_{\partial B_r}|\nabla G|^{n-2}\frac{\partial G}
{\partial n} =1-\dint_{B_r}G^{n-1}dx=1+O(r^{n}\log^{n-1}r).$$ Then,
we get (\ref{3.11}).
\par \medskip
We have
\begin{equation}\label{q}\dint_{\mathbb{R}^n\setminus B_\delta}u_k^\frac{n}{n-1}\Phi'(\beta_ku_k^\frac{n}{n-1})dx
\leq C\dint_{\mathbb{R}^n\setminus B_\delta}u_k^ndx \rightarrow 0\ .
\end{equation}
Recall that $U_k\in H^{1,n}_0(B_{R_k})$. By equation (\ref{3.10}) we
get
$$\dint_{\mathbb{R}^n\setminus B_\delta}(|\nabla U_k|^n+U_k^n)dx=
\frac{c_k^\frac{n}{n-1}}{\lambda_k}\dint_{\mathbb{R}^n\setminus
B_\delta}
{u_k^\frac{n}{n-1}}\Phi'(\beta_ku_k^\frac{n}{n-1})dx-\dint_{\partial
B_\delta} \frac{\partial U_k}{\partial n}|\nabla U_k|^{n-2}U_kdS.$$
By (\ref{q}) and (\ref{3.8a}) we then get
$$\begin{array}{ll}\
     \lim\limits_{k\rightarrow+\infty}\dint_{\mathbb{R}^n\setminus B_\delta}
       (|\nabla U_k|^n+U_k^n)dx
       &=-\lim\limits_{k\rightarrow+\infty}\dint_{\partial B_\delta}
            \frac{\partial U_k}{\partial n}|\nabla U_k|^{n-2}U_kdS \vspace{0.2cm}\\[\mv]
       &=-G(\delta)\dint_{\partial B_\delta}
          \frac{\partial G}{\partial n}|\nabla G|^{n-2}dS \vspace{0.2cm}\\[\mv]
       &=G(\delta)(1-\dint_{B_\delta}G^{n-1}dx).
  \end{array}$$

$\hfill\Box$\\

\par\bigskip
We are now in the position to complete the proof of Theorem 1.1: We
have seen in (\ref{3.8b}) that
$$\dint_{\mathbb{R}^n\setminus B_R}\Phi(\beta_ku_k^\frac{n}{n-1})dx
\leq C.$$ So, we only need to prove on $B_R$,
$$\dint_{B_R}e^{\beta_ku_k^\frac{n}{n-1}}dx<C$$

The classical Trudinger-Moser inequality implies that
$$\dint_{B_R}e^{\beta_k((u_k-u_k(R))^+)^\frac{n}{n-1}}dx<C=C(R).$$
By Proposition \ref{t3.7}, $u_k(R)=O(\frac{1}{c_k^\frac{1}{n-1}})$,
and hence we have
$$u_k^\frac{n}{n-1}\leq ((u_k-u_k(R))^++u_k(R))^\frac{n}{n-1}
\leq ((u_k-u_k(R))^+)^\frac{n}{n-1}+C_1,$$ Then, we get
$$\dint_{B_R}e^{\beta_ku_k^\frac{n}{n-1}}\leq C\,'.$$

$\hfill\Box$\\

\section{The proof of  Proposition \ref{t1.3}}

We will use a result of Carleson and Chang (see \cite{C-C}):

\begin{lem}\label{t4.1}Let $B$ be the unit ball in $\mathbb{R}^n$.
Assume that $u_k$ is a sequence in  $H^{1,n}_0(B)$ with
$\int_B|\nabla u_k|^ndx=1$. If $u_k\rightharpoondown 0$, then
$$\limsup_{k\rightarrow+\infty}\dint_B(e^{\alpha_n|u_k|^\frac{n}{n-1}}-1)dx
\leq|B|e^{1+1/2+\cdots+1/(n-1)}.$$
\end{lem}

{\it Proof of Proposition {\ref{t1.3}}:} Set
$u_k'(x)=\frac{(u_k(x)-u_k(\delta))^+}{\|\nabla
u_k\|_{L^n(B_\delta)}}$ which is in $H^{1,n}_0(B_\delta)$. Then by
the result of Carleson and Chang, we have
$$\limsup_{k\rightarrow+\infty}\dint_{B_\delta}e^{\beta_k{u_k'}^\frac{n}{n-1}}
\leq|B_\delta|(1+e^{1+1/2+\cdots+1/(n-1)}).$$

By Lemma \ref{t3.8}, we have
$$\dint_{\mathbb{R}^n\setminus B_\delta}
(|\nabla
c_k^\frac{1}{n-1}u_k|^n+(c_k^\frac{1}{n-1}u_k)^n)dx\rightarrow
 G(\delta)(1-\dint_{B_\delta}
G^{n-1}dx)\ ,
$$
and therefore we get
\begin{equation}\label{n}
\dint_{B_\delta}|\nabla u_k|^ndx=1-\dint_{\mathbb{R}^n\setminus
B_\delta} (|\nabla u_k|^n+u_k^n)dx-\dint_{B_\delta}u_k^ndx
=1-\frac{G(\delta) +\epsilon_k(\delta)}{c_k^\frac{n}{n-1}}\, ,
\end{equation}
where $\lim\limits_{\delta\rightarrow 0}\lim\limits_{k\rightarrow
+\infty} \epsilon_k(\delta)=0$.
\par\medskip
By (\ref{3.8b}) in Lemma \ref{t3.5} we have
$$\lim_{L\rightarrow+\infty}\lim_{k\rightarrow+\infty}\dint_{B_\rho\setminus B_{Lr_k}}
e^{\beta_ku_k^\frac{n}{n-1}}dx=|B_\rho|\, ,
$$
for any $\rho<\delta$. Furthermore, on $B_\rho$ we have by (\ref{n})
$$\begin{array}{lll}
  (u_k')^\frac{n}{n-1}\le
  \displaystyle \frac{u_k^\frac{n}{n-1}}
      {(1-\frac{G(\delta)+\epsilon_k(\delta)}{c_k^\frac{n}{n-1}})^{\frac{1}{n-1}}}
    &=& \displaystyle
    u_k^\frac{n}{n-1}(1 + \frac{1}{n-1}\frac{G(\delta)+\epsilon_k(\delta)}{c_k^\frac{n}{n-1}}+O(\frac{1}
         {c_k^\frac{2n}{n-1}}))\vspace{0.2cm}\\[\mv]
    &=& \displaystyle
    u_k^\frac{n}{n-1} + \frac{1}{n-1}G(\delta)(\frac{u_k}{c_k})^\frac{n}{n-1}
       +O(c_k^\frac{-n}{n-1})\vspace{0.3cm}\\[\mv]
    &\leq& \displaystyle
    u_k^\frac{n}{n-1}-\frac{\log\delta^n}{(n-1)\alpha_n} \ .
  \end{array}$$\\
Then we have
$$\lim_{L\rightarrow+\infty}\lim_{k\rightarrow+\infty}\dint_{B_\rho\setminus B_{Lr_k}}
e^{\beta_k{u_k'}^\frac{n}{n-1}}dx\leq
O(\delta^{-n})\lim_{L\rightarrow+\infty}
\lim_{k\rightarrow+\infty}\dint_{B_\rho\setminus B_{Lr_k}}
e^{\beta_ku_k^\frac{n}{n-1}}dx\rightarrow |B_\rho|O(\delta^{-n}).$$
Since $u_k'\rightarrow 0$ on $B_\delta\setminus B_\rho$, we get
$$\lim_{k\rightarrow+\infty}
\dint_{B_\delta\setminus
B_\rho}(e^{\beta_k{u_k'}^\frac{n}{n-1}}-1)dx =0,$$ then
$$0\leq\lim_{L\rightarrow+\infty}\lim_{k\rightarrow+\infty}
\dint_{B_\delta\setminus
B_{Lr_k}}(e^{\beta_k{u_k'}^\frac{n}{n-1}}-1)dx\leq
|B_\rho|O(\delta^{-n}).$$ Letting $\rho\rightarrow 0$, we get
$$\lim_{L\rightarrow+\infty}\lim_{k\rightarrow+\infty}
\dint_{B_\delta\setminus
B_{Lr_k}}(e^{\beta_k{u_k'}^\frac{n}{n-1}}-1)dx =0.$$ So, we have
$$\lim_{L\rightarrow+\infty}\lim_{k\rightarrow+\infty}
\dint_{B_{Lr_k}}(e^{\beta_k{u_k'}^\frac{n}{n-1}}-1)dx\leq
e^{1+1/2+\cdots+1/(n-1)} |B_\delta|.$$

Now, we fix an $L$. Then for any $x\in B_{Lr_k}$, we have
$$\begin{array}{lll}
   \beta_ku_k^\frac{n}{n-1}&=&
        \beta_k(\frac{u_k}{\|\nabla u_k\|_{L^n(B_\delta)}})^\frac{n}{n-1}(\dint_{B_\delta}|\nabla u_k|^n
        dx)^\frac{1}{n-1}\\[\mv]
     &=&\beta_k(u_k'+\frac{u_k(\delta)}{\|\nabla u_k\|_{L^n(B_\delta)}})^\frac{n}{n-1}
       (\dint_{B_\delta}|\nabla u_k|^ndx )^\frac{1}{n-1} \vspace{0.2cm}\\[\mv]
         &&\s\ \big( \hbox{using that }u_k(\delta)=O(\frac{1}{c_k^\frac{1}{n-1}})
   \hbox{ and }\|\nabla u_k\|_{L^n(B_\delta)}=1+O(\frac{1}{c_k^\frac{n}{n-1}}) \, \big) \vspace{0.2cm}\\[\mv]
     &=&\beta_k\big( u_k'+u_k(\delta)+O(\frac{1}{c_k^\frac{n+1}{n-1}})\big)^\frac{n}{n-1}\
        \big( \dint_{B_\delta}|\nabla u_k|^ndx \big)^\frac{1}{n-1}\vspace{0.2cm}\\[\mv]
     &=&\beta_k {u_k'}^{\frac{n}{n-1}}\big( 1+\frac{u_k(\delta)}{u_k'}+O(\frac{1}{c_k^\frac{2n}{n-1}})\big)^\frac{n}{n-1}\
     (1-\frac{G(\delta)+\epsilon_k(\delta)}{c_k^\frac{n}{n-1}}
        )^\frac{1}{n-1}\\[\mv]
     &=&\beta_k{u_k'}^{\frac{n}{n-1}}\left[1+\frac{n}{n-1}\frac{u_k(\delta)}{u_k'}
         -\frac{1}{n-1}\frac{G(\delta)+\epsilon_k(\delta)}{c_k^\frac{n}{n-1}}+O(\frac{1}{c_k^\frac{2n}{n-1}})\right]\ .
\end{array}$$
It is easy to check that
$$\frac{u_k'(r_k x)}{c_k}\rightarrow 1\ , \s \hbox{ and } \s
{\big( u_k'(r_kx)\big)}^\frac{1}{n-1}u_k(\delta)\rightarrow
G(\delta) \ .
$$
So, we get
$$\begin{array}{ll}
    \lim\limits_{L\rightarrow+\infty}\lim\limits_{k\rightarrow+\infty}
    \dint_{B_{Lr_k}}(e^{\beta_ku_k^\frac{n}{n-1}}-1)dx&=
           \lim\limits_{L\rightarrow+\infty}\lim\limits_{k\rightarrow+\infty}
           e^{\alpha_nG(\delta)}
           \dint_{B_{Lr_k}}(e^{\beta_k{u_k'}^\frac{n}{n-1}}-1)dx \vspace{0.3cm}\\[\mv]
          & \displaystyle \leq
           e^{\alpha_nG(\delta)}\ \delta^n\frac{\omega_{n-1}}{n}\
           e^{1+1/2+\cdots+1/(n-1)}\vspace{0.3cm} \\[\mv]
           & \displaystyle =
           e^{\alpha_n(-\frac1{\alpha_n}\log\delta^n + A + O(\delta^n\log^n\delta))}\
           \delta^n\frac{\omega_{n-1}}{n}\
           e^{1+1/2+\cdots+1/(n-1)}.
  \end{array}$$
\par \medskip \noindent
Letting $\delta\rightarrow 0$, then the above inequality together
with Lemma \ref{t3.2} imply  Proposition \ref{t1.3}.
\par \bigskip

\section{The test function 1}

In this section, we will  construct a function sequence
$\{u_\epsilon\}\subset H^{1,n}(\mathbb{R}^n)$ with
$\|u_\epsilon\|_{H^{1,n}}=1$ which satisfies
$$\dint_{\mathbb{R}^n}\Phi(\alpha_n|u_\epsilon|^\frac{n}{n-1})dx>
\frac{\omega_{n-1}}{n}e^{\alpha_nA+1+1/2+
 \cdots+/1(n-1)},$$
 for $\epsilon>0$ sufficiently small.

Let
$$u_\epsilon=\left\{\begin{array}{ll}
                 C-\frac{(n-1)\log(1+c_n|\frac{x}
                    {\epsilon}|^\frac{n}{n-1})+{\Lambda_\epsilon}}{\alpha_n
                   {C}^\frac{1}{n-1}}&|x|\leq L\epsilon\\[\mv]
                 \frac{G(|x|)}{{C}^\frac{1}{n-1}}&
                   |x|> L\epsilon\ ,
               \end{array}\right.
$$
\par \medskip \noindent
where ${\Lambda_\epsilon}, {C} $ and ${L}$ are functions of
$\epsilon$ (which will be defined later, by \big(\ref{5.1}),
(\ref{5.2}), (\ref{5.5})\big) which satisfy
\par \medskip

i) \ ${L}\rightarrow +\infty$, ${C}\rightarrow +\infty$, and
${L}\epsilon\rightarrow 0$, as $\epsilon\rightarrow 0$\ ;

ii) \
${C}-\frac{(n-1)\log(1+c_n{L}^\frac{n}{n-1})+{\Lambda_\epsilon}}{\alpha_n{C}^\frac{1}{n-1}}
=\frac{G(L\epsilon)}{{C}^\frac{1}{n-1}}$\ ;

iii)\ $\frac{\log{L}}{C^\frac{n}{n-1}}\rightarrow 0$, as
$\epsilon\rightarrow 0$\ .\\

We use the normalization of $u_\epsilon$ to obtain information on
$\Lambda_\epsilon, C $ and $L$. We have
$$
\begin{array}{lll}
     \dint_{\mathbb{R}^n\setminus B_{L\epsilon}}(|\nabla u_\epsilon|^n+u^n_\epsilon)dx&=&
        \frac{1}{{C}^\frac{n}{n-1}}\ \big(\dint_{B_{L\epsilon}^c}
      |\nabla G|^ndx+
       \int_{B_{L\epsilon}^c}G^ndx\big)\vspace{0.2cm}\\[\mv]
     &=&\frac{1}{{C}^\frac{n}{n-1}}
        \dint_{\partial B_{L\epsilon}}G(L\epsilon)
     |\nabla {G}|^{n-2}\frac{\partial G}
       {\partial n}dS\vspace{0.2cm}\\[\mv]
     &=&\displaystyle \frac{G(L\epsilon)-G(L\epsilon)\dint_{B_{L\epsilon}}
         G^{n-1}dx}{{C}^\frac{n}{n-1}}\ .
  \end{array}$$
and
$$\begin{array}{llll}
     \dint_{B_{L\epsilon}}|\nabla u_\epsilon|^ndx&=& \displaystyle\frac{n-1}{\alpha_n{C}^\frac{n}{n-1}}
        \dint_0^{c_n{L}^\frac{n}{n-1}} \frac{u^{n-1}}{(1+u)^n}du\vspace{0.3cm} \\[\mv]
      &=&\displaystyle \frac{n-1}{\alpha_n{C}^\frac{n}{n-1}}\dint_{0}^{c_n{L}^\frac{n}{n-1}}
         \frac{((1+u)-1)^{n-1}}{(1+u)^n}du \vspace{0.3cm} \\[\mv]
      &=& \displaystyle \frac{n-1}{\alpha_n{C}^\frac{n}{n-1}}\sum_{k=0}^{n-2}\frac{C_{n-1}^k(-1)^{
          n-1-k}}{n-k-1}\vspace{0.2cm}\\[\mv]
         && \displaystyle  +\frac{n-1}{\alpha_n{C}^\frac{n}{n-1}}\log(1+c_nL^\frac{n}{n-1})
            +O(\frac{1}{L^\frac{n}{n-1}{C}^\frac{n}{n-1}}) \vspace{0.3cm}\\[\mv]
      &=&\displaystyle  -\frac{n-1}{\alpha_n{C}^\frac{n}{n-1}}\big( 1+1/2+1/3+\cdots+1/(n-1)\big)\vspace{0.2cm}\\[\mv]
         &&\displaystyle  +\frac{n-1}{\alpha_n{C}^\frac{n}{n-1}}\log(1+
         c_n{L}^\frac{n}{n-1})+O(\frac{1}{{L}^\frac{n}{n-1}
         {C}^\frac{n}{n-1}})\ ,
 \end{array}$$
where we used the fact
$$
-\sum_{k=0}^{n-2}\frac{C_{n-1}^k(-1)^{
          n-1-k}}{n-k-1}=1+\frac{1}{2}+\cdots+\frac{1}{n-1}\ .$$
It is easy to check that
$$
\dint_{B_{L\epsilon}}|u_\epsilon|^ndx=O((L\epsilon)^nC^n\log L)\ ,
$$
and thus we get
$$\begin{array}{lll}
    \dint_{\mathbb{R}^n}(|\nabla u_\epsilon|^n+u_\epsilon^n)dx
    &=& \displaystyle \frac{1}{\alpha_n {C}^\frac{n}{n-1}}\Big\{-(n-1)\big(1+1/2+\cdots+1/(n-1)\big) +
    \alpha_n A \vspace{0.2cm}\\[\mv]
       && \displaystyle +(n-1)\log(1+c_n{L}^\frac{n}{n-1})
       -\log({L}\epsilon)^n+\phi \Big\}\ ,
   \end{array}$$
where
$$\phi=O\left((L\epsilon)^nC^n\log L+(L\epsilon)^n\log^n{L\epsilon}+L^\frac{-n}
{n-1}\right).$$

Setting $\int_{\mathbb{R}^n}(|\nabla
u_\epsilon|^n+u_\epsilon^n)dx=1$, we obtain
\begin{equation}\label{5.1}
  \begin{array}{lll}
    \alpha_n{C}^\frac{n}{n-1}&=&-(n-1)\big(1+1/2+\cdots+1/(n-1)\big)+\alpha_nA
     +\log\frac{(1+c_n{L}^\frac{n}{n-1})^{n-1}}{{L}^n}-
     \log{\epsilon^n}+\phi\\[\mv]
    &=&-(n-1)\big(1+1/2+\cdots+1/(n-1)\big)+\alpha_nA+\log{\frac{\omega_{n-1}}{n}}-
      \log{\epsilon^n}+\phi\ .
  \end{array}
\end{equation}
By ii) we have
$${\alpha_n{C}^\frac{n}{n-1}}-(n-1)\log(1+c_n{L}^\frac{n}{n-1})+{\Lambda_\epsilon}
= \alpha {G(L\epsilon)}
$$
and hence
$$
-(n-1)\big(1+1/2+\cdots+1/(n-1)\big)+\alpha_nA
     -\log{(L\epsilon)^n}+\phi +{\Lambda_\epsilon}
= \alpha {G(L\epsilon)} \ ;
$$
this implies that
\begin{equation}\label{5.2}
{\Lambda_\epsilon}=-(n-1)(1+1/2+\cdots+1/(n-1))+\phi\ .
\end{equation}
\par \bigskip
Next, we compute $
\dint_{B_{L\epsilon}}e^{\alpha_n|u_\epsilon|^\frac{n}{n-1}}dx$\ .
\par \medskip \noindent
Clearly, $\varphi(t)=|1-t|^\frac{n}{n-1}+\frac{n}{n-1}t$ \ is
increasing when $0\leq t\leq 1$ and decreasing when $t\leq 0$, then
$$
|1-t|^\frac{n}{n-1}\geq 1-\frac{n}{n-1}t\ ,\ \hbox{ when } \ |t|<1\
.
$$
Thus we have by ii), for any $x\in B_{L\epsilon}$
\begin{equation}\label{5.3}
  \begin{array}{lll}
    \alpha_nu_\epsilon^\frac{n}{n-1}
    &=&\displaystyle \alpha_n{C}^\frac{n}{n-1}\Big|1-\frac{(n-1)\log(1+c_n|\frac{x}{\epsilon}|^\frac{n}{n-1})
      +{\Lambda_\epsilon}}
     {\alpha_n{C}^\frac{n}{n-1}}\Big|^\frac{n}{n-1}\\[\mv]
    &\geq&\displaystyle \alpha_n{C}^\frac{n}{n-1}(1-\frac{n}{n-1}
     \frac{(n-1)\log(1+c_n|\frac{x}{\epsilon}|^\frac{n}{n-1})+{\Lambda_\epsilon}}
     {\alpha_n{C}^\frac{n}{n-1}}).
  \end{array}
\end{equation}

Then we have
$$\begin{array}{lll}
    \displaystyle \dint_{B_{L\epsilon}}e^{\alpha_n|u_\epsilon|^\frac{n}{n-1}}dx
    &\geq&\displaystyle \dint_{B_{L\epsilon}} e^{\alpha_n{C}^\frac{n}{n-1}
     -n\log(1+c_n|\frac{x}{\epsilon}|^\frac{n}{n-1})-\frac{n}{n-1}{\Lambda_\epsilon}}
     \vspace{0.2cm}\\[\mv]
    &=&\displaystyle e^{\alpha_n{C}^\frac{n}{n-1}
     -\frac{n}{n-1}{\Lambda_\epsilon}}\int_{B_L}\frac{\epsilon^n}{(1+c_n|x|^\frac{n}
     {n-1})^n}dx \vspace{0.2cm} \\[\mv]
    &=& \displaystyle e^{\alpha_n{C}^\frac{n}{n-1}
     -\frac{n}{n-1}{\Lambda_\epsilon}}(n-1)\epsilon^n\int_{0}^{c_n{L}^\frac{n}{n-1}}
     \frac{u^{n-2}}{(1+u)^n}du \vspace{0.2cm}\\[\mv]
    &=& \displaystyle e^{\alpha_n{C}^\frac{n}{n-1}
     -\frac{n}{n-1}{\Lambda_\epsilon}}(n-1)\epsilon^n\int_{0}^{c_n{L}^\frac{n}{n-1}}
     \frac{((u+1)-1)^{n-2}}{(1+u)^n}du \vspace{0.2cm}\\[\mv]
    &=& \displaystyle e^{\alpha_n{C}^\frac{n}{n-1}
     -\frac{n}{n-1}{\Lambda_\epsilon}}\epsilon^n(1+O({L}^{-\frac{n}{n-1}})) \vspace{0.2cm}\\[\mv]
    &=&\displaystyle \frac{\omega_{n-1}}{n}e^{\alpha_nA+1+1/2+\cdots+1/(n-1)} \vspace{0.2cm}\\[\mv]
    &&+\ O\left((L\epsilon)^nC^n\log L+L^\frac{-n}{n-1}+(L\epsilon)^n\log^n{L\epsilon}\right).
  \end{array}$$
Here, we used the fact
$$\sum_{k=0}^m\frac{(-1)^{m-k}}{m-k+1}C_m^k=\frac{1}{m+1}\ .
$$
Then
$$
\dint_{B_{L\epsilon}}\Phi(\alpha_nu_\epsilon^\frac{n}{n-1})dx\geq\frac{\omega_{n-1}}{n}
e^{\alpha_nA+1+1/2+\cdots+1/(n-1)}+
  O\left((L\epsilon)^nC^n\log L+L^\frac{-n}{n-1}+(L\epsilon)^n\log^n{L\epsilon}\right)\ .
$$
\par \medskip \noindent
Moreover, on $ \mathbb{R}^n\setminus B_{L\epsilon}$ we have the
estimate
$$
\dint_{\mathbb{R}^n\setminus
B_{L\epsilon}}\Phi(\alpha_nu_\epsilon^\frac{n}{n-1})dx \geq
\frac{\alpha_n^{n-1}}{(n-1)!}\dint_{\mathbb{R}^n\setminus
B_{L\epsilon}} \left|\frac{G(x)}{C^\frac{1}{n-1}}\right|^ndx\ ,
$$
and thus we get
\begin{equation}\label{5.4}
  \begin{array}{l}
    \dint_{\hspace{-0.2cm}\mathbb{R}^n}\Phi(\alpha_nu_\epsilon^\frac{n}{n-1})dx
       \geq\frac{\omega_{n-1}}{n}e^{\alpha_nA+1+1/2+\cdots+1/(n-1)}\\[\mv]
      \s + \ \frac{\alpha_n^{n-1}}{(n-1)!}\dint_{\mathbb{R}^n\setminus B_{L\epsilon}}
             \left|\frac{G(x)}{C^\frac{1}{n-1}}\right|^ndx+
         O\left((L\epsilon)^nC^n\log L+L^\frac{-n}{n-1}+(L\epsilon)^n\log^n{L\epsilon}
         \right)\vspace{0.5cm}\\[\mv]
      \s =\ \displaystyle \frac{\omega_{n-1}}{n}e^{\alpha_nA+1+1/2+\cdots+1/(n-1)}\\[\mv]
      \s + \ \frac{\alpha_n^{n-1}}{(n-1)!\; C^\frac{n}{n-1}}
         \left[\dint_{\mathbb{R}^n\setminus B_{L\epsilon}}\hspace{-0.2cm}|G(x)|^ndx +
         O\Big((L\epsilon)^nC^{n+\frac{n}{n-1}}\log L+\frac{C^\frac{n}{n-1}}{L^\frac{n}{n-1}}
          + C^\frac{n}{n-1}(L\epsilon)^n\log^n{L\epsilon}\Big)\right]
 \end{array}
\end{equation}
\par \bigskip

We now set
\begin{equation}\label{5.5}
{L}=-\log{\epsilon}\ ;
\end{equation}
then $L\epsilon\rightarrow 0$ as $\epsilon\rightarrow 0$. We then
need to prove that there exists a $C = C(\epsilon)$ which solves
equation (\ref{5.1}). We set
$$
f(t)=-\alpha_nt^\frac{n}{n-1}-(n-1)(1+1/2+\cdots+1/(n-1))+\alpha_nA+\log{\frac{\omega_{n-1}}
{n}}-\log{\epsilon^n}+\phi\ ,
$$
Since
$$f((-\frac{2}{\alpha_n}\log\epsilon^n)^\frac{n}{n-1})=\log\epsilon^n+o(1)+\phi<0$$
for $\epsilon$ small, and
$$f((-\frac{1}{2\alpha_n}\log\epsilon^n)^\frac{n}{n-1})=-\frac{1}{2}
\log\epsilon^n+o(1)+\phi>0$$ for $\epsilon$ small, $f$ has a zero in
$\left((-\frac{1}{2\alpha_n} \log\epsilon^n)^\frac{n-1}{n},
(-\frac{2}{\alpha_n}\log\epsilon^n)^\frac{n-1}{n}\right)$. Thus, we
defined $C$, and it satisfies
$$\alpha_nC^\frac{n}{n-1}=-\log\epsilon^n+O(1).$$
Therefore, as $\epsilon\rightarrow 0$, we have
$$\frac{\log L}{C^\frac{n}{n-1}}\rightarrow 0,$$
and then
$$(L\epsilon)^nC^{n+\frac{n}{n-1}}\log L+C^\frac{n}{n-1}L^\frac{-n}{n-1}+
C^\frac{n}{n-1}(L\epsilon)^n\log^n{L\epsilon}\rightarrow 0.$$
Therefore, i), ii), iii) hold and we can conclude from (\ref{5.4})
that for $\epsilon > 0$ sufficiently small
$$\int_{\mathbb{R}^n}\Phi(\alpha_nu_\epsilon^\frac{n}{n-1})dx>\frac{\omega_{n-1}}{n}e^{
  \alpha_nA+1+1/2+\cdots+1/(n-1)}\ .
$$
\par \bigskip

\section{The test function 2}
In this section we construct, for $n > 2$, functions $u_\epsilon$
such that
$$
\dint_{\mathbb{R}^n}\Phi(\alpha_n(\frac{u_\epsilon}{\|u_\epsilon\|_{H^{1,n}}})^\frac{n}
{n-1})dx> \frac{\alpha_n^{n-1}}{(n-1)!}\ ,
$$
for $\epsilon > 0$ sufficiently small.
\par \bigskip
Let  $\epsilon^n=e^{-\alpha_nc^\frac{n}{n-1}}$, and
$$
u_\epsilon=\left\{\begin{array}{ll}
             c\quad & \quad |x|<L\epsilon\\[\mv]
             \frac{-n\log\frac{x}{L}}{\alpha_nc^\frac{1}{n-1}}\quad & \quad L\epsilon\leq|x|\leq L\\[\mv]
             0 \quad & \quad L\leq |x| \ ,
           \end{array}\right.$$
where $L$ is a function of $\epsilon$ which will be defined later.
\par \bigskip
We have
$$\dint_{\mathbb{R}^n}|\nabla u_\epsilon|^n=1,$$
and
$$\dint_{\mathbb{R}^n}u_\epsilon^ndx=\frac{\omega_{n-1}}{n}c^n(L\epsilon)^n+
       \frac{\omega_{n-1}n^nL^n}{\alpha_n^nc^\frac{n}{n-1}}\dint_{\epsilon}^1r^{n-1}\log^nrdr.$$

Then
$$\begin{array}{lll}
    \dint_{\mathbb{R}^n}\Phi(\alpha_n(\frac{u_\epsilon}{\|u_\epsilon\|_{H^{1,n}}})^\frac{n}{n-1})dx
        &\geq& \displaystyle \frac{\alpha_n^{n-1}}{(n-1)!}\
        \frac{\int_{\mathbb{R}^n}u_\epsilon^ndx}{1 + \int_{\mathbb{R}^n}u_\epsilon^ndx} +
           \frac{\alpha_n^n}{n!}\ \frac{\int_{\mathbb{R}^n\setminus B_{L\epsilon}}
            u_\epsilon^\frac{n^2}{n-1}}{(1+\int_{\mathbb{R}^n}u_\epsilon^ndx)^\frac{n}{n-1}}
           dx \vspace{0.3cm}\\[\mv]
        &=&\displaystyle \frac{\alpha_n^{n-1}}{(n-1)!}-\frac{\alpha_n^{n-1}}{(n-1)!}\frac{1}{1+\frac{
          \omega_{n-1}}{n}c^n(L\epsilon)^n+\frac{\omega_{n-1}n^nL^n}{\alpha_n^nc^\frac{n}{n-1}}
          \int_{\epsilon}^1
          r^{n-1}\log^nrdr} \vspace{0.3cm} \\[\mv]
        && \displaystyle + \frac{\alpha_n^n}{n!}\frac{\omega_{n-1}L^n/c^\frac{n^2}{(n-1)^2}(\frac{n}{\alpha_n}
        )^\frac{n^2}{n-1}
          \int_\epsilon^1r^{n-1}\log^\frac{n^2}{n-1}r}{\Big(1+\frac{
          \omega_{n-1}}{n}c^n(L\epsilon)^n+\frac{\omega_{n-1}n^nL^n}{\alpha_n^nc^\frac{n}{n-1}
          }\int_{\epsilon}^1
          r^{n-1}\log^nrdr \Big)^\frac{n}{n-1}} \vspace{0.3cm} \\[\mv]
  \end{array}
$$

We now ask that $L$ satisfies
$$\frac{c^\frac{n}{n-1}}{L^n}\rightarrow 0,\s as\s \epsilon\rightarrow 0.\eqno (6.1)$$

Then, for sufficiently small $\epsilon$, we have
$$\begin{array}{l}
   \displaystyle -\frac{\alpha_n^{n-1}}{(n-1)!}\frac{1}{1+\frac{
          \omega_{n-1}}{n}c^n(L\epsilon)^n+\frac{\omega_{n-1}n^nL^n}{\alpha_n^nc^\frac{n}{n-1}
          }\int_{\epsilon}^1 r^{n-1}\log^nrdr}\ + \vspace{0.2cm} \\
        + \ \displaystyle \frac{\alpha^n}{n!}\frac{\omega_{n-1}L^n/c^\frac{n^2}
        {(n-1)^2}(\frac{n}{\alpha_n})^\frac{n^2}{n-1}
          \int_\epsilon^1r^{n-1}\log^\frac{n^2}{n-1}r}{\Big(1+\frac{\omega_{n-1}}{n}c^n(L\epsilon)^n +
          \frac{\omega_{n-1}n^nL^n}{\alpha_n^nc^\frac{n}{n-1}}\int_{\epsilon}^1
          r^{n-1}\log^nrdr \Big)^\frac{n}{n-1}}\vspace{0.2cm}\\[\mv]
     \s\s\s \geq  \displaystyle B_1L^{n-\frac{n^2}{n-1}}-B_2\frac{c^\frac{n}{n-1}}{L^n}\vspace{0.2cm} \\
     \s\s\s =\ \displaystyle \frac{
      c^\frac{n}{n-1}}{L^n}(B_1\frac{L^{2n -\frac{n^2}{n-1}}}{c^\frac{n}{n-1}}-B_2) \vspace{0.2cm} \\
     \s\s\s =  \ \displaystyle \frac{
      c^\frac{n}{n-1}}{L^n}(B_1 \frac{L^{\frac n{n-1}(n-2)}}{c^{\frac n{n-1}}} - B_2) \ ,
  \end{array}$$
where $B_1$, $B_2$ are positive constants.

When $n>2$, we may choose $L= b\, c^{\frac 1 {n-2}}$; then, for $b$
sufficiently large, we have
$$
B_1 \frac{L^{\frac n{n-1}(n-2)}}{c^{\frac n{n-1}}} - B_2 = B_1\,
b^{\frac n{n-1}(n-2)} - B_2 > 0 \ ,
$$
and (6.1) holds. Thus, we have proved that for $\epsilon > 0$
sufficiently small
$$
\dint_{\mathbb{R}^n}\Phi(\alpha_n(\frac{u_\epsilon}{\|u_\epsilon\|_{H^{1,n}(\mathbb{R}^n)}})^\frac{n}{n-1})dx
 > \frac{\alpha_n^{n-1}}{(n-1)!}\ .
$$

\par 
\noindent Yuxiang  Li\\
{\it ICTP, Mathematics Section, Strada Costiera 11, I-34014 Trieste,
Italy}\\
{\it E-mail address: liy@ictp.it}\\

\noindent Bernhard Ruf \\
{\it Dipartimento di Matematica, Universit$\grave{a}$ di Milano, via
Saldini 50, 20133 Milan, Italy}\\ {\it E-mail address:
Bernhard.Ruf@mat.unimi.it}

\end{document}